\documentclass[
11pt, article]%
{revtex4}

\begin{document}

\title{Inverse problem for determining the order of the fractional derivative in equations of mixed type}

\author{{R.~R.}~{Ashurov}}
\affiliation{Institute of Mathematics of the Academy of  Sciences, Republic of Uzbekistan}

\author{R.~T.{Zunnunov}}
\affiliation{Institute of Mathematics of the Academy of  Sciences, Republic of Uzbekistan}

\begin{abstract} 
In this paper  the inverse problem of determining the fractional
orders in mixed-type  equations is considered. In one part of the
domain, the considered equation is the subdiffusion equation with
a fractional derivative in the sense of Gerasimov-Caputo of the
order $\alpha\in (0,1)$ and in the other part - a wave equation
with a fractional derivative of the order  $\beta\in (1,2)$. The
elliptic part of the equation is a second-order operator,
considered in   $N$-dimensional domain $\Omega$. Assuming the
parameters $\alpha$ and $\beta$ to be unknown, additional
conditions are found that provide an unambiguous determination of
the required parameters.

\end{abstract}

\keywords{mixed-type  equation, fractional derivative in the sense of Gerasimov-Caputo, the inverse problem of determining the orders  of  fractional derivative.} 

\maketitle

\section{Introduction}

The theory of differential equations with fractional derivatives has gained considerable popularity and importance in the last few decades, mainly due to its applications in numerous seemingly distant fields of science and technology (see, for example,\cite{{as1}, {as2}, {as3}, {as4}, {as5}, {as6}, {as7}, {as8}}). The theory of equations of mixed-type with fractional derivatives is much less developed (see \cite{as9} and the references in it).

Inverse problems in the theory of partial differential equations
(of integer and fractional order) are such problems in which
together with the solution of the differential equation, it is
required to determine a coefficient(s) of the equation, or/and the
right-hand side (source function). Interest in the study of
inverse problems is due to the importance of their applications in
various branches of mechanics, seismology, medical tomography, and
geophysics (see, the monographs of M.M. Lavrent'ev et al.
\cite{as10}, S.I. Kabanikhin \cite{{as11},{as12}}). In
references \cite {{as13}, {as14}, {as15}, {as16}, {as17}, {as18}} various inverse problems for
equations of integer order were considered, and in \cite{{as19}, {as20}, {as21}, {as22}} - inverse problems to determine the right-hand side
for equations of fractional order were considered.

In considering of fractional-order equation as a model equation for analyzing various processes, the order of the fractional derivative is often unknown and difficult to measure directly. Inverse problems for determining this unknown parameter are of theoretical interest and are necessary for solving initial-boundary value problems and studying the properties of solutions.
This is a relatively new type of inverse problem, which arises only when considering equations of fractional order (see the survey article \cite{as23} and the references in it, and the recent articles \cite{{as24}, {as25}, {as26}, {as27}, {as28}}).

In this article, we investigate the inverse problem of determining the orders of fractional derivatives in the sense of Gerasimov-Caputo in equations of mixed-type. As far as we know, similar inverse problems for equations of mixed-type are studied for the first time. However, it should be noted that initial-boundary value problems and inverse problems for determining the right-hand side in equations of mixed-type with the fractional derivatives were investigated in \cite{as9}, \cite{as29}, \cite{as30}. In these studies, ordinary differential expressions were considered as the elliptic part and the existence of a weak solution was studied.

Let us move on to exact descriptions of the objects of study.

Let $\Omega$ be an arbitrary $N$-dimensional domain with a sufficiently smooth boundary $\partial\Omega$. Let further the second-order differential operator
\[
A(x,D)u=\sum\limits_{i,j=1}^N D_i[a_{i,j}(x)D_ju]-c(x)u
\]
 be a symmetric elliptic operator in $\Omega$, i.e
\[
a_{i,j} (x) = a_{j,i}(x)\quad \text{and} \quad \sum_{i,j=1}^N
a_{i,j}(x) \xi_i \xi_j\geq a\sum_{i,j=1}^N  \xi^2_i,
\]
for all $x\in \Omega$ and  $\xi_i$, where $a=const>0$ and $D_ju=\frac{\partial u}{\partial x_j}, \ j=1, \dots, N.$
 Consider the spectral problem with the Dirichlet condition

\begin{equation}\label{ash1}
-A(x,D)v(x) =\lambda v(x), \quad x\in \Omega;
\end{equation}
\begin{equation}\label{ash2}
v(x) = 0, \quad x\in \partial \Omega.
\end{equation}
 It is known (see, \cite{as31}, p.100) that if the coefficients of operator $A(x, D)$ and the boundary of the domain $\Omega$ are sufficiently smooth and  $c(x)> 0$, then the spectral problem (\ref{ash1}) - (\ref{ash2}) has a complete in $L_2(\Omega)$ set of  orthonormal eigen functions       $\{v_k(x)\}$, $k\geq 1$   and a countable set of positive eigenvalues $\{\lambda_k\}$.

 To determine the fractional part of our differential equation we first introduce a fractional integral in the sense of Riemann-Liouville of order $\alpha<0$ of the function $f$ defined on $[0, \infty)$, according to the following formula (see, \cite{as7}, p.14)
$$
\partial_t^\alpha f(t)=\frac{1}{\Gamma
    (-\alpha)}\int\limits_0^t\frac{f(\xi)}{(t-\xi)^{\alpha+1}} d\xi,
\quad t>0,
$$
 provided that the right-hand side of equality exists. Here $\Gamma(t)$ is Euler's gamma function. Using the fractional integral, we define the fractional derivative in the sense of Gerasimov-Caputo of order a $\alpha$, $k-1<\alpha\leq k$, $k\in \mathbf{N}$ by setting (see \cite{as7}, p.14 and \cite{as32})
 \begin{equation}\label{ash3}
 D_{0t}^\alpha f(t)= \partial_t^{\alpha-k}\frac{d^k}{dt^k} f(t).
 \end{equation}

  If $\alpha=k$, then the fractional derivative coincides with the ordinary classical derivative: $D_{0t}^k f(t) =\frac{d^k}{dt^k} f(t)$.

 Let $0<\alpha < 1$ and $1<\beta<2$. In the domain  $\Omega \times (-T,
 +\infty)$, $T>0$, consider an equation of mixed-type
\begin{equation}\label{ash4}
\left\{
\begin{array}{c}
D_{0t}^\alpha u(x,t) -A(x,D) u(x,t) = 0, \quad x\in \Omega, \quad
t>0;\\
D^\beta_{t0} u(x,t) -A(x,D) u(x,t) = 0, \quad x\in \Omega, \quad
-T<t<0.
\end{array}
\right.
\end{equation}
We take the Dirichlet condition as the boundary condition, i.e
\begin{equation}\label{ash5}
u(x,t)=0, \quad x\in \partial\Omega, \quad t\geq -T.
\end{equation}
 Let the initial condition be of the following form
\begin{equation}\label{ash6}
u(x, -T)=\varphi(x),\quad x\in \overline{\Omega}.
\end{equation}
We take the gluing conditions in the form
\begin{equation}\label{ash7}
u(x, +0)=u(x, -0), \quad \lim\limits_{t\rightarrow +0}
D^\alpha_{0t}u(x, t)= u_t(x, -0), \quad x\in \Omega.
\end{equation}
If  $\alpha =1$, these last conditions mean the continuity of the
solution  $u(x, t)$ and its derivative with respect to  $t$ on the
line of variation of the type of equation:  $t=0$.

The initial-boundary value problem (\ref{ash4})-(\ref{ash7}) is called  \emph{the forward problem}.

Further, it is convenient to introduce the notation: $G^+=\overline{\Omega} \times (0,
+\infty)$ and $G^-=\overline{\Omega} \times (-T, 0)$.

\textbf{Definition} Function $u(x, t)\in C(\overline{\Omega}\times
    [-T,0)\cap(0, \infty))$ with properties
    \begin{enumerate}
        \item$A(x,D) u(x,t)\in
        C(G^+\bigcup G^-)$,
        \item
        $D_{0t}^\alpha u(x,t)\in C(G^+)$
        \item$D_{t0}^\beta u(x, t)\in C(G^-) $
    \end{enumerate}
 satisfying conditions (\ref{ash4})-(\ref{ash7}) is called the (classical) solution of the forward problem.

Note that the requirement of continuity in the closed domain $\overline{\Omega}$ of all derivatives of the solution  $u(x,t)$ included in equation (\ref{ash4}) is not caused by the essence. However, on the one hand, the uniqueness of just such a solution can be proved quite simply, and on the other hand, the solution determined by the Fourier method satisfies the conditions indicated above.

Let $E_{\rho, \mu}(t)$ be the Mittag-Leffler function
$$
E_{\rho, \mu}(t)= \sum\limits_{n=0}^\infty \frac{t^n}{\Gamma(\rho
    n+\mu)}.
$$
The following function is denoted by  $\Delta_k$, $k\geq 1$:
$$
\Delta_k\equiv \Delta_k(T,\beta)=\lambda_k T
E_{\beta,2}\big(-\lambda_k T^\beta\big)- E_{\beta,1}\big(-\lambda_k T^\beta\big).
$$

\textbf{Lemma}There exists a constant  $T_0= T_0(\lambda_1, \beta)$ such that for  $T\geq
    T_0$ the following estimate is true
    \[
    \Delta_k>\delta_0 >0, \quad k\geq 1,
    \]
where the constant  $\delta_0= \delta_0(T)$ does not depend on $\lambda_k$.

\

We denote the classical Sobolev space by $W^m_2(\Omega)$ and the integer part of $b$ by $[b]$. The solution of the forward problem in domains $G^\pm$ is denoted by $u^\pm(x,t)$, and by $\varphi_k$- the Fourier coefficients of the function $\varphi (x)$ by the system of eigenfunctions $\{v_k(x)\}$ of the spectral problem (\ref{ash1})-(\ref{ash2}).

\textbf{Theorem}\label{ft} Let the condition of Lemma \ref{delta} be satisfied and let $\varphi$ satisfy the following conditions
    \begin{equation}\label{ash8}
    \varphi(x)\in W_2^{\big[\frac{N}{2}\big]+2}(\Omega),
    \end{equation}
    \begin{equation}\label{ash9}
    \varphi(x)=A(x,D)\varphi(x)= \cdot\cdot\cdot\,
    =A^{\big[\frac{N}{4}\big]}(x,D)\varphi(x)=0, \quad x\in \partial
    \Omega.
    \end{equation}
    Then there is a unique solution to the forward problem and it can be represented in the form of the following series
        \begin{equation}\label{ash10}
    u^+(x,t)=\sum\limits_{k=1}^\infty \frac{ E_{\alpha, 1} \big(-\lambda_k
        t^\alpha\big)\varphi_k v_k(x)}{\Delta_k},\quad 0\leq t<\infty,
    \end{equation}
    \begin{equation}\label{ash11}
    u^-(x,t)=\sum\limits_{k=1}^\infty \frac{ \bigg[|t|
        \lambda_kE_{\beta,2} \big(-\lambda_k |t|^\beta\big)-E_{\beta, 1}
        \big(-\lambda_k |t|^\beta\big)\bigg]\varphi_k v_k(x)}{\Delta_k}, \quad
    -T\leq t\leq 0,
    \end{equation}

    which converge absolutely and uniformly in $x\in
    \overline{\Omega}$ and in $t$ in the domains indicated above.

\ 

Now assume that the orders of fractional derivatives $\alpha$ and $\beta$ are unknown and consider the inverse problem of determining these parameters. In this case, it is natural to assume that $\varphi(x)$ is not a null function, since otherwise, by virtue of Theorem  \ref{ft}, the solution of the forward problem is $u(x,t)\equiv 0$ for any $\alpha$ and $\beta$, which does not allow determining these parameters. Since there are two unknowns, we set the following two additional conditions:

\begin{equation}\label{ash12}
\int\limits_\Omega |u(x,t_1)|^2 dx = d_1,
\end{equation}
\begin{equation}\label{ash13}
\int\limits_\Omega u(x,-t_2) v_{k_0}(x) dx=d_2,
\end{equation}
 where $d_1, d_2$ - are the numbers given arbitrarily, $t_1, t_2$ are the positive numbers, defined in
 Theorem 4 given below, and $k_0\geq 1$ is an arbitrary integer such that $\varphi_{k_0} \neq 0$
 (obviously, such numbers exist since $\varphi(x)$ is not identically zero).

The initial-boundary value problem (\ref{ash4})-(\ref{ash7}) with additional conditions (\ref{ash12}), (\ref{ash13}) is called the \emph{inverse problem} for determining parameters $\alpha$ and $\beta$. If $u(x,t)$ is a solution to the forward problem and parameters $\alpha$ and $\beta$ satisfy conditions (\ref{ash12}), (\ref{ash13}) then the triple   $\{u(x,t),\alpha,\beta\}$ is called a \emph{ solution of the inverse problem} .

When solving the inverse problem, we assume that
$$
0<\alpha_1\leq \alpha<1, \quad 1<\beta_1\leq \beta\leq \beta_2<2.
$$
where $\alpha_1$, $\beta_1$ and $\beta_2$  are the numbers given
arbitrarily from the corresponding intervals. It should be noted
that the integral in (\ref{ash12}) depends on both parameters
$\alpha$ and $\beta$, while the left-hand side of condition
(\ref{ash13}) depends only on parameter $\beta$ (see the solution
to the forward problem, (\ref{ash10}) and  (\ref{ash11})).
Therefore, to solve the inverse problem, we first determine
parameter $\beta^\star$ from condition (\ref{ash13}), then,
assuming that parameter $\beta^\star$ is already known, we
determine parameter $\alpha^\star$ satisfying condition
(\ref{ash12}).

Let us introduce the following notation
\[
W(\alpha, \beta)=\int\limits_\Omega |u(x,t_1)|^2 dx, \quad
P(\beta)=  \frac{\Delta_{k_0}(t_2, \beta)}{\Delta_{k_0}(T,
    \beta)}.
\]
Using the explicit form of solution (\ref{ash11}) and the orthonormality of functions  $v_k(x)$, condition (\ref{ash13}) can be rewritten in the following form (note, that $\varphi_{k_0}\neq 0$)
\begin{equation}\label{ash14}
P(\beta) \varphi_{k_0}=d_2.
\end{equation}

For solvability of the inverse problem, it is obvious that numbers $d_1$ and $d_2$ cannot be specified in an arbitrary way. A necessary condition for the solvability of equation (\ref{ash13}) is the fulfillment of the following inequality (in the last section it was proven that function $P(\beta)$ is monotonically increasing, and function $W(\alpha, \beta^\star)$ is monotonically decreasing)

\begin{equation}\label{ash15}
P(\beta_1) \leq \frac{d_2}{\varphi_{k_0}} \leq P(\beta_2).
\end{equation}
 In other words, if the reverse inequality holds, then under no value of $\beta\in [\beta_1, \beta_2]$ there is a solution $u(x,t)$ of the forward problem that satisfies condition (\ref{ash13}). In what follows, we assume that this condition is satisfied. The condition relative to number $d_1$ under which the inverse problem has a solution is given in the following theorem.

\textbf{Theorem 2.}\label{it}  Let the conditions of Theorem \ref{ft} be satisfied. There is a number $T_{1,2}= T_{1,2}(\lambda_1, \beta_1,
    \beta_2)$ such that for $T_{1,2} \leq t_2<T $  there is a unique number of $\beta^\star$, satisfying condition (\ref{ash13}). Further, there exists a positive number $T_3=T_3(\lambda_1, \alpha_1)$ such that for $t_1\geq T_3$ the necessary and sufficient condition for the existence of a unique solution $\{u(x,t),\alpha,\beta^\star\}$ of the inverse problem has the following form
    \begin{equation}\label{ash16}
    W(1, \beta^\star)< d_1\leq W(\alpha_1, \beta^\star).
    \end{equation}


 {\bf Remark 1.} It is seen that condition (\ref{ash12}) is given in the form of the norm of the solution in $L_2(\Omega)$, and condition (\ref{ash13}) - is in the form of the Fourier coefficient of the solution to the forward problem. As follows from the proof of Theorem \ref{it}  condition (\ref{ash12}) can also be specified in the form of the Fourier coefficient, and condition (\ref{ash13}) - in the form of the norm in $L_2(\Omega)$.

In the next section, we prove the uniqueness of the solution to the forward problem. Section 3 is devoted to the study of the question of the existence of a solution to this problem. Applying the Fourier method, the solution is presented in the form of a formal series in the eigenfunctions of the spectral problem (\ref{ash1}) - (\ref{ash2}). Moreover, the conditions on function $\varphi$ that guarantee the convergence and differentiability of these series are first given in terms of belonging to the domain of fractional powers of self-adjoint extensions of operator $A(x, D)$. Then, using the fundamental result of V.A. Il'in  \cite{as31}  , we formulate these conditions in terms of Sobolev spaces. The fourth section is devoted to the study of the monotonicity property of the Mittag-Leffler functions  $E_{\rho, \mu} (t)$ with respect to parameter $\rho$. In the last section, the inverse problem is considered and Theorem \ref{it} is proven.

\section{Uniqueness}

Let us first prove Lemma  \ref{delta}, that is, we establish a lower estimate for the quantity
$$
\Delta_k\equiv \Delta_k(T, \beta)=\lambda_k T
E_{\beta,2}\big(-\lambda_k T^\beta\big)- E_{\beta,1}\big(-\lambda_k T^\beta\big).
$$
For this, we use the asymptotic estimates of the Mittag - Leffler
functions ( see \cite{as6}, p.134, and Lemmas \ref{e2} and
\ref{e1} below) for large values of the arguments. Let
$1<\beta<2$. Then for $\mu\geq 1$ the following estimates hold:
\[
\mu E_{\beta,2}(-\mu)\ =\ \frac{1}{\Gamma(2-\beta)}\ +\
\frac{r_1(\mu, \beta)}{\mu},
\]
\[
E_{\beta,1}(-\mu)\ =\ \frac{r_2(\mu, \beta)}{\mu},
\]
where
\[
|r_1(\mu, \beta)|\ \leq\ C_1, \,\, \text{and}\,\,|r_2(\mu, \beta)|\
\leq\ C_2,
\]
constants $C_1$ and $C_2$  depend only on  $\beta$.

Therefore, for  $T\geq1$ we have
\[
\frac{\mu E_{\beta,2}(-\mu)}{T^{\beta-1}}\ -\ E_{\beta,1}(-\mu)\ =
\]
\begin{equation}\label{ash17}
=\ \frac{1}{T^{\beta-1}\Gamma(2-\beta)}\ +\ \frac{1}{\mu}
\left[\frac{r_1(\mu, \beta)}{T^{\beta-1}}\ -\ r_2(\mu,
\beta)\right].
\end{equation}

\textbf{Lemma 2.} Let $1<\beta<2$ and $T\geq1$. Assume that
    $$
    \mu_0(T)\ =\ \max\left\{1,\ 2\,\Gamma(2-\beta) \left(C_1\ +\ C_2
    T^{\beta-1}\right)\right\}.
    $$
    Then for  $\mu\geq \mu_0$ the following estimate is true:
    \begin{equation}\label{ash18}
    \frac{\mu E_{\beta,2}(-\mu)}{T^{\beta-1}}\ -\ E_{\beta,1}(-\mu)\
    \geq \ \delta_0(T),
    \end{equation}
    where
    $$
    \delta_0(T)\ =\  \frac{1}{2\,T^{\beta-1}\Gamma(2-\beta)}.
    $$

\textbf{Proof.} In virtue of equation  (\ref{ash17}), it is sufficient to show the existence of such number
    $\mu_0\geq1$, that for $\mu\geq \mu_0$ the following estimate is true
    \[
    \frac{1}{\mu} \left|\frac{r_1(\mu, \beta)}{T^{\beta-1}}\ -\
    r_2(\mu, \beta)\right|\ \leq\ \frac{1}{2\,
        T^{\beta-1}\Gamma(2-\beta)}\ .
    \]

    Using the estimates for  $r_1$ and  $r_2$, we have
    \[
    \frac{1}{\mu} \left|\frac{r_1(\mu, \beta)}{T^{\beta-1}}\ -\
    r_2(\mu, \beta)\right|\ \leq\
    \frac1\mu\left[\frac{C_1}{T^{\beta-1}}\ +\ C_2 \right].
    \]
    Therefore it is sufficient to show that
    \[
    \frac1\mu\left[\frac{C_1}{T^{\beta-1}}\ +\ C_2 \right] \leq\
    \frac{1}{2\, T^{\beta-1}\Gamma(2-\beta)}.
    \]
    But this estimate is obviously true for  $\mu\geq \mu_0(T)$.
Lemma 1 easily follows from this lemma. Indeed, let
\[
\mu_k\ =\ \lambda_k T^{\beta}, \quad k=1, 2, ...
\]

Since $\beta >1$, then there exists  $T_0 = T_0(\lambda_1,\beta)$
such that for  $T\geq T_0$ we have  $\lambda_1 T^\beta\ \geq\
\mu_0(T)$. Then for these  $T$
\[
\mu_k\ \geq\ \mu_1\ =\ \lambda_1 T^\beta\ \geq\ \mu_0(T)
\]
for all  $k\geq 1$. Consequently, estimate (\ref{ash18}) implies
\[
\frac{\mu_k E_{\beta,2}(-\mu_k)}{T^{\beta-1}}\ -\
E_{\beta,1}(-\mu_k)\ \geq \ \delta_0(T),
\]
which, by the definition of  $\mu_k$, coincides with the assertion of Lemma
\ref{delta}.

We turn to the proof of the uniqueness of the solution to the forward problem. Suppose that under the conditions of Theorem  \ref{ft} there are two solutions: $u_1(x,t)$ and $u_2(x,t)$. Let us prove
that $u(x,t)=u_1(x,t)-u_2(x,t)\equiv 0$. Since the problem under consideration is linear, to determine  $u(x,t)$ we have the initial-boundary value problem  (\ref{ash4}) - (\ref{ash7}), and the initial condition has the following form:
\begin{equation}\label{ash19}
u(x, -T)=0, \quad x\in \overline{\Omega}.
\end{equation}

Let $u(x,t)$ satisfy all the conditions of this problem and let $v_k$ is arbitrary eigenfunction of the spectral problem
(\ref{ash1})-(\ref{ash2}) with the corresponding eigenvalue $\lambda_k$. Consider the following function
\begin{equation}\label{ash20}
w_k(t)=\int\limits_\Omega u(x,t)v_k(x)dx.
\end{equation}
From equations  (\ref{ash4}) we obtain (by virtue of Definition \ref{def} differentiation is possible under the integral sign)
$$
D_{0t}^\alpha w_k(t)=\int\limits_\Omega D_{0t}^\alpha
u(x,t)v_k(x)dx= \int\limits_\Omega A(x, D)u(x,t)v_k(x)dx, \quad
t>0,
$$
$$
D_{t0}^\beta w_k(t)=\int\limits_\Omega D_{t0}^\beta
u(x,t)v_k(x)dx= \int\limits_\Omega A(x, D)u(x,t)v_k(x)dx, \quad
-T<t<0.
$$
Taking into account boundary conditions  (\ref{ash2}) and (\ref{ash5}), the integration is done by parts. Then, by virtue of equation (\ref{ash1}), we have
\begin{equation}\label{ash21}
D_{0t}^\alpha w_k(t)+ \lambda_k w_k(t)=0, \quad t>0,
\end{equation}
\begin{equation}\label{ash22}
D_{t0}^\beta w_k(t)+\lambda_k w_k(t)=0, \quad -T<t<0.
\end{equation}
Let us rewrite the gluing conditions  (\ref{ash7}) for function  $w_k(t)$  in the following form:
$$
w_k(+0) =w_k(-0),\quad \lim\limits_{t\rightarrow +0} D_{0t}^\alpha
w_k(t) =w'_k(-0).
$$
First, assume that  $w_k(+0)$, $w_k(-0)$ and $w'_k(-0)$ are known given numbers. Then the solutions of the Cauchy problems related to equations  (\ref{ash21}) and (\ref{ash22}), have the following form, respectively ( see \cite{as7} p.16, \cite{as33})
\begin{equation}\label{ash23}
\left\{
\begin{array}{c}
 w_k(t)= w_k(+0) E_{\alpha, 1}\big(-\lambda_k t^\alpha\big), \quad t>0,\\
w_k(t)= w_k(-0) E_{\beta, 1}\big(-\lambda_k |t|^\beta\big)+ w'_k(-0)|t|
E_{\beta, 2}\big(-\lambda_k |t|^\beta\big), \quad -T<t<0.
\end{array}
\right.
\end{equation}

 From the explicit form of function  $w_k(t)$ for $t>0$ it follows that
$\lim\limits_{t\rightarrow +0} D_{0t}^\alpha w_k(t) = -\lambda_k
w_k(+0).$ Hence, by virtue of the gluing conditions and initial condition
(\ref{ash19}), we obtain a system of two equations
to determine two numbers $w_k(-0)$ and $w'_k(-0)$
(note that the third unknown number  $w_k(+0)$ is determined from the following condition  $w_k(+0)=w_k(-0)$) we obtain a system of two equations
\begin{equation}\label{ash24}
\left\{
\begin{array}{c}
\lambda_k w_k(-0) +w'_k(-0)=0;\\
E_{\beta,1}\big(-\lambda_k T^\beta\big) w_k(-0) + T
E_{\beta,2}\big(-\lambda_k T^\beta\big) w'_k(-0)=0.
\end{array}
\right.
\end{equation}
Since the determinant of this system is  $\Delta_k> \delta_0>0$,
the system has a unique solution:
$w_k(+0)=w_k(-0)=w'_k(-0)=0$. This, in turn, means that
$w_k(t)\equiv 0$, $t\geq -T$, for all  $k\geq 1$. Therefore, due to the completeness of the system of eigenfunctions  $\{v_k(x)\}$, we have $u(x,t) = 0$ for all $x\in \overline{\Omega}$ and $t>-T$. So, the uniqueness of the solution to the forward problem is proven.

{\bf Remark 2.}It should be noted that if $ T $ is not chosen
large enough, then the direct problem either has more than one
solution or has no solution at all. Indeed, it follows from the
definition of the Mittag-Leffler function that for small $
\lambda_k T $ the function $ \Delta_k (T, \beta) $ is negative.
Then, by Lemma \ref{delta}, some of the $ \Delta_k (T, \beta) $
vanish for the corresponding $ T $. And if the determinant of the
system (\ref{ash24}) vanishes, then for such $ T $ the Fourier
coefficients $ w_k (t) $ of the function $ u (x, t) $ in the
system $ \{v_k (x) \} $ are either not unique or do not exist
depending on the augmented matrix.

\section{EXISTENCE}

Let us prove the existence of a solution to the forward problem by the classical Fourier method. For this, we first construct a formal solution, assuming that the initial function $\varphi(x)$ expands in a Fourier series.

According to the Fourier method, function $u(x,t)$ is sought in the form of a formal series
\begin{equation}\label{ash25}
u(x,t) = \sum\limits_{k=1}^\infty w_k(t) v_k(x),
\end{equation}
here functions  $w_k(t)$ are the solutions of equations  (\ref{ash21}) and
(\ref{ash22}) and can be represented in the form of  (\ref{ash23}). Using initial condition  (\ref{ash6}) and the gluing condition (\ref{ash7}), for the unknown numbers  $w_k(-0)$ and $w'_k(-0)$ ( note that the third unknown number  $w_k(+0)$ is determined from condition
$w_k(+0)=w_k(-0)$) we obtain a system of two equations
\begin{equation}\label{ash26}
\left\{
\begin{array}{c}
\lambda_k w_k(-0) +w'_k(-0)=0;\\
E_{\beta,1}\big(-\lambda_k T^\beta\big) w_k(-0) + T
E_{\beta,2}\big(-\lambda_k T^\beta\big) w'_k(-0)=\varphi_k.
\end{array}
\right.
\end{equation}
For sufficiently large $T$, according to Lemma \ref{delta}, the determinant of system (\ref{ash26}) is strictly positive. Solving this system, we obtain representations (\ref{ash10}) and (\ref{ash11}) for the formal solution (\ref{ash25}). It remains to substantiate the Fourier method. For this, we need to introduce some definitions.

Let $\tau$ - be an arbitrary real number. In the space
$L_2(\Omega)$ we introduce the operator  $\hat{A}^\tau$, acting according to the following rule
$$
\hat{A}^\tau g(x)= \sum\limits_{k=1}^\infty \lambda_k^\tau g_k
v_k(x), \quad g_k=(g,v_k).
$$
Obviously, the given operator  $\hat{A}^\tau$ with the domain
$$
D(\hat{A}^\tau)=\bigg\{g\in L_2(\Omega):  \sum\limits_{k=1}^\infty
\lambda_k^{2\tau} |g_k|^2 < \infty\bigg\}
$$
is a self-adjoint operator. If by  $A$ we denote an operator in
$L_2(\Omega)$ acting according to the rule  $Ag(x)=A(x,D) g (x)$
and with domain  $D(A)=\bigg\{g\in C^2(\overline{\Omega}) :g(x)=0, x\in
\partial\Omega\bigg\}$, then the operator $\hat{A}\equiv\hat{A}^1$ is a
self-adjoint extension in $L_2(\Omega)$ of operator  $A$.

Our reasoning further will be largely based on the methods
developed by M.A. Krasnoselsky at al. \cite{as34}. The key role in
this technique is played  the following lemma, proven in the
same book  (\cite{as34}, p.453):

\textbf{Lemma.}\label{Kr} Let a multi-index $\sigma$ be such that   $|\sigma|\leq 2$ and $q > \frac{|\sigma|}{2}+\frac{N}{4}$. Then the operator  $D^\sigma \hat{A}^{-q}$  acts (completely) continuously from  $L_2(\Omega)$ to $C(\overline{\Omega})$ and the following estimate is true
    $$
    ||D^\sigma \hat{A}^{-q} g||_{C(\Omega)} \leq C
    ||g||_{L_2(\Omega)}.
    $$


We assume that the Fourier coefficients of function $\varphi$ satisfy the following condition
\begin{equation}\label{ash27}
\sum\limits_{k=1}^\infty \lambda_k^{2(\tau+1)} |\varphi_k|^2 \leq
C_\varphi<\infty
\end{equation}
for some  $\tau> \frac{N}{4}$ and show that functions
(\ref{ash10}) and (\ref{ash11}) satisfy conditions  1-3 of Definition \ref{def}.

Consider the partial sum of series  (\ref{ash11})

\[
u_j^-(x,t)=\sum\limits_{k=0}^j \frac{ \Big[|t|
    \lambda_kE_{\beta,2} \big(-\lambda_k |t|^\beta\big)-E_{\beta, 1}
    \big(-\lambda_k |t|^\beta\big)\Big]\varphi_k v_k(x)}{\Delta_k}, \quad
-T\leq t< 0.
\]
Since $\hat{A}^{-\tau-1} v_k(x) = \lambda_k^{-\tau-1} v_k(x)$,
we have
\[
u_j^-(x,t)=\hat{A}^{-\tau-1}\sum\limits_{k=0}^j \frac{ \Big[|t|
    \lambda_kE_{\beta,2} \big(-\lambda_k |t|^\beta\big)-E_{\beta, 1}
    \big(-\lambda_k |t|^\beta\big)\Big]\lambda_k^{\tau+1}\varphi_k
    v_k(x)}{\Delta_k}.
\]
Applying Lemma \ref{Kr}, we obtain

$$
||D^\sigma u^-_j||_{C(\Omega)}=\bigg|\bigg|D^\sigma
\hat{A}^{-\tau-1}\sum\limits_{k=0}^j \frac{ \Big[|t|
    \lambda_kE_{\beta,2} \big(-\lambda_k |t|^\beta\big)-E_{\beta, 1}
    \big(-\lambda_k |t|^\beta\big)\Big]\lambda_k^{\tau+1}\varphi_k
    v_k(x)}{\Delta_k}\bigg|\bigg|_{C(\Omega)}
$$
\begin{equation}\label{ash28}
\leq C \bigg|\bigg|\sum\limits_{k=0}^j \frac{ \Big[|t|
    \lambda_kE_{\beta,2} \big(-\lambda_k |t|^\beta\big)-E_{\beta, 1}
    \big(-\lambda_k |t|^\beta\big)\Big]\lambda_k^{\tau+1}\varphi_k
    v_k(x)}{\Delta_k}\bigg|\bigg|_{L_2(\Omega)}.
\end{equation}
Using the orthonormality of system $\{v_k\}$ and Lemma \ref{delta}, we have
$$
||D^\sigma u^-_j||^2_{C(\Omega)}\leq \frac{C}{\delta_0^2}
\sum\limits_{k=0}^j \Big|\Big[|t| \lambda_kE_{\beta,2} \big(-\lambda_k
|t|^\beta\big)-E_{\beta, 1} \big(-\lambda_k
|t|^\beta\big)\Big]\lambda_k^{\tau+1}\varphi_k \Big|^2.
$$
Further, taking into account the estimate of the Mittag-Leffler function (see, for example \cite{as7}, p. 13)
\begin{equation}\label{ash29}
|E_{\rho, \mu}(-t)|\leq \frac{C}{1+t}, \quad t>0,
\end{equation}
where $\mu$ is an arbitrary complex number; we obtain
$$
||D^\sigma u^-_j||^2_{C(\Omega)}\leq \frac{C}{\delta_0^2}
\sum\limits_{k=1}^j\bigg|\frac{1+|t|\lambda_k}{1+|t|^\beta\lambda_k}\bigg|^2
\lambda_k^{2(\tau+1)}|\varphi_k|^2\leq C
C_\varphi\delta_0^{-2}\big(1+t^{-2(\beta-1)}\big).
$$

This implies the uniform in $x\in\overline{\Omega}$ convergence for all $t\in[-T, 0)$  of the sum (\ref{ash11}) differentiated in variables  $x_i$. On the other hand, sum (\ref{ash28}) converges for any permutation of its terms since these terms are mutually orthogonal. This implies the absolute convergence of the differentiated series (\ref{ash11}) in variables $x_i$.

Repeating completely similar reasoning, we can state that the same assertions are true for series (\ref{ash10}), which determines function $u^+(x,t)$.

Now let us verify the validity of the term-by-term application of operators $D_{0t}^\alpha$ and $D_{t0}^\beta$  to series (\ref{ash10}) and (\ref{ash11}), respectively. Consider the second of these operators since the first is considered similarly (see (\ref{ash31})). It is easy to check that

$$
D_{t0}^\beta\sum\limits_{k=1}^j w_k(t)v_k(x)=
-\sum\limits_{k=1}^j\lambda_k
w_k(t)v_k(x)=-A(x,D)\hat{A}^{-\tau-1}\sum\limits_{k=1}^j
\lambda^{\tau+1}_k w_k(t)v_k(x), \quad -T< t<0.
$$
Absolute and uniform convergence of the last sum has already been proven above.

It remains to check the fulfillment of the gluing conditions (\ref{ash7}). The fact that these conditions are formally satisfied follows from the choice of functions $w_k(x)$ (see (\ref{ash26})). Therefore, it is sufficient to show the continuity in domain  $\Omega$ of functions $u(x,t)$ at $t=0$, $D_{0t}^\alpha u^+(x,t)$  at $t=+0$ and $u^-_x(x,t)$ at  $t=-0$.

Applying Lemma \ref{Kr}  for $\sigma =0$ and repeating the above reasoning (see (\ref{ash28})), we obtain

$$
||u^+_j||_{C(\Omega)}=\bigg|\bigg|\hat{A}^{-\tau}\sum\limits_{k=0}^j
\frac{ E_{\alpha, 1} \big(-\lambda_k
    t^\alpha\big)\lambda_k^{\tau}\varphi_k
    v_k(x)}{\Delta_k}\bigg|\bigg|_{C(\Omega)}
$$
\begin{equation}\label{ash30}
\leq C \bigg|\bigg|\sum\limits_{k=0}^j \frac{ E_{\alpha, 1}
    \big(-\lambda_k t^\alpha\big)\lambda_k^{\tau}\varphi_k
    v_k(x)}{\Delta_k}\bigg|\bigg|_{L_2(\Omega)}\leq
\frac{C}{\delta_0^2}\sum\limits_{k=0}^j \bigg|E_{\alpha, 1}
\big(-\lambda_k t^\alpha\big)\lambda_k^{\tau}\varphi_k\bigg|^2.
\end{equation}
 By virtue of the estimate for the Mittag - Leffler function (\ref{ash29}), we finally have
\[
||u^+_j||_{C(\Omega)}\leq C C_\varphi\delta_0^{-2},\,\, t\geq 0,
\]
which means the uniform convergence in domain  $(x,t)\in
\overline{\Omega}\times [0,\infty)$ of the series (\ref{ash10}), consisting of continuous functions. The continuity of function $u^-(x,t)$ for $t\rightarrow -0$ is proven in a similar way.

Since $D_{0t}^\alpha E_{\alpha, 1} (-\lambda_k t^\alpha)
=-\lambda_k E_{\alpha, 1} (-\lambda_k t^\alpha)$
(see \cite{as3}, formula (4.3.1)), similar to (\ref{ash30}) we have

$$
||D_{0t}^\alpha
u^+_j||_{C(\Omega)}=\bigg|\bigg|\hat{A}^{-\tau}\sum\limits_{k=0}^j
\frac{ -E_{\alpha, 1} \big(-\lambda_k
    t^\alpha\big)\lambda_k^{\tau+1}\varphi_k
    v_k(x)}{\Delta_k}\bigg|\bigg|_{C(\Omega)}\leq
$$
\begin{equation}\label{ash31}
\leq \frac{C}{\delta_0^2}\sum\limits_{k=0}^j \bigg|E_{\alpha, 1}
\big(-\lambda_k t^\alpha\big)\lambda_k^{\tau+1}\varphi_k\bigg|^2\leq C
C_\varphi\delta_0^{-2},\,\, t\geq 0,
\end{equation}
which means uniform convergence in domain $(x,t)\in
\overline{\Omega}\times [0,\infty)$ of the differentiated series (\ref{ash10}). This, in particular, implies the continuity of function $D_{0t}^\alpha u^+(x,t)$ for $t\rightarrow + 0$.

 We turn to the proof of the continuity of derivative $u^-_x(x,t)$. We have the following formula (see  \cite{as3}, formula (4.3.1))
\[
\frac{d}{dt} \big[t \lambda_kE_{\beta,2} \big(-\lambda_k
t^\beta\big)-E_{\beta, 1} \big(-\lambda_k t^\beta\big)\big]=\lambda_k
E_{\beta,1} \big(-\lambda_k t^\beta\big)-\frac{1}{t} E_{\beta,0}
\big(-\lambda_k t^\beta\big),\,\, 0<t<T.
\]
On the other hand, it is easy to see that
\[
\frac{1}{t} E_{\beta,0} \big(-\lambda_k t^\beta\big)=-t^{\beta-1}\lambda_k
E_{\beta,\beta} \big(-\lambda_k t^\beta\big).
\]
Therefore, applying Lemma \ref{Kr} for $\sigma =0$ and the estimate for the Mittag - Leffler function (\ref{ash29}), we finally have
$$
||u^-_x(x,-t)||_{C(\Omega)}=\bigg|\bigg|\hat{A}^{-\tau}\sum\limits_{k=0}^j
\frac{ \Big[E_{\beta,1} \big(-\lambda_k t^\beta\big)+t^{\beta-1}
    E_{\beta,\beta} \big(-\lambda_k
    t^\beta\big)\Big]\lambda_k^{\tau+1}\varphi_k
    v_k(x)}{\Delta_k}\bigg|\bigg|_{C(\Omega)}
$$
\[ \leq \frac{C}{\delta_0^2}\sum\limits_{k=0}^j
\bigg|\Big[E_{\beta,1} \big(-\lambda_k t^\beta\big)+t^{\beta-1}
E_{\beta,\beta} \big(-\lambda_k
t^\beta\big)\Big]\lambda_k^{\tau+1}\varphi_k\bigg|^2\leq C
C_\varphi\delta_0^{-2}T^{2(\beta-1)},\,\, 0\leq t\leq T.
\]

Thus, we have shown that under condition (\ref{ash27}) all the assertions of Theorem \ref{ft} hold. On the other hand, in the fundamental study by V. A. Il'in \cite{as31}, it was proven that conditions (\ref{ash8}) and (\ref{ash9}) of Theorem \ref{ft} guarantee the convergence of the series (\ref{ash27}).

Thus, the proof of Theorem \ref{ft} is complete.

\section{PROPERTIES OF THE MITTAG-LEFFLER FUNCTIONS }

In this section, we study the behavior of the Mittag-Leffler functions $E_{\alpha, 1}(-\lambda t^\alpha)$, $E_{\beta,
    1}(-\lambda t^\beta)$  and  $E_{\beta, 2}(-\lambda t^\beta)$ depending on parameters  $\alpha$ and $\beta$  respectively. These functions participate in series (\ref{ash10}) and (\ref{ash11}), which determine the solution of the forward problem. When evaluating these functions, various constants arise, depending on some parameter $b$. All these constants are denoted by one symbol $C(b)$.

\textbf{Lemma.}\label{e2}
    Let $\beta_1\geq 1$ and $\beta_2<2$. There exists a positive number  $T_2=T_2(\lambda_1, \beta_2)$ such that for all  $t\geq
    T_2$ and $\lambda \geq \lambda_1$ function $e_{\lambda,2}(\beta) =
    E_{\beta, 2}(-\lambda t^\beta)$ is positive and monotonically decreasing in  $\beta\in [\beta_1, \beta_2]$ and the following estimates are true
    \[
    e_{\lambda,2}(\beta) =\frac{1}{\lambda
        t^\beta}\cdot\frac{(2-\beta)}{\Gamma
        (3-\beta)}+\frac{O(1)}{(\lambda t^\beta)^2},
    \]
    \[
    \frac{d}{d\beta} e_{\lambda,2}(\beta) =- \frac{1}{\lambda
        t^\beta}\cdot\frac{(2-\beta)[\ln t - \Psi (3-\beta)]-1}{ \Gamma
        (3-\beta)}+\frac{O(1)}{(\lambda t^\beta)^2},
    \]
where  $O(1)$ depends only on $\beta_2$.


\textbf{Proof. }
    Let us denote by  $ \delta (\theta) $ a contour traversed in the direction of non-decreasing  $ \arg \zeta $ and consisting of the following parts: ray $ \arg \zeta = - \theta $ , $ | \zeta | \geq
    1 $, arc $ (- \theta) \leq \arg \zeta \leq \theta $, $ | \zeta | = 1
    $ and ray $ \arg \zeta = \theta $, $ | \zeta | \geq 1 $. If $ 0
    <\theta <\pi $, then the contour  $  \delta (\theta) $ divides the entire complex  $ \zeta $ - plane into two unbounded parts, namely  $ G^{(-)} (\theta) $ to the left of  $ \delta (\theta) $ and $
    G^{(+)} (\theta) $ to the right of it. Contour  $ \delta (\theta) $
    is called the Hankel path (see \cite{as6}, p.126).

    Parameter $\theta$ is taken in the following form
    \[
    \theta=\theta_\epsilon \beta, \quad \theta_\epsilon =\frac{\pi}{2}
    +\epsilon, \quad \epsilon =\frac{\pi -\frac{\pi}{2}
        \beta_2}{1+\beta_2}>0, \quad \beta\in [\beta_1, \beta_2].
    \]
    Then   $\frac{\pi}{2}\beta<\theta\leq \pi -\epsilon<\pi$ and
    $-\lambda t^\beta\in G^{(-)} (\theta) $ and by the definition of contour $\delta(\theta)$, we have (see  \cite{as6}, p.135)
    \begin{equation}\label{ash32}
    E_{\beta, 2}\big(-\lambda t^\beta\big)= \frac{1}{\lambda t^\beta \Gamma
        (2-\beta)}-\frac{1}{2\pi i \beta\lambda
        t^\beta}\int\limits_{\delta(\theta)}\frac{e^{\zeta^{1/\beta}}\zeta^{1-1/\beta}}{\zeta+\lambda
        t^\beta} d\zeta = f_1(\beta, t)+f_2(\beta, t).
    \end{equation}

Let us first prove that derivative  $\frac{d}{d\beta}
e_{\lambda,2} (\beta)$ is negative for all $\beta\in [\beta_1, \beta_2]$. From now on we denote the derivative in $\beta$ of  the function $f(\beta, t)$ through  $f'(\beta, t)$.

Function $f'_1(\beta, t)$ is easily calculated. Indeed, let $\Psi(\beta)$ be the logarithmic derivative of function $\Gamma(\beta)$ (the definition and properties of function $\Psi$  see in \cite{as35}). Then $\Gamma'(\beta) = \Gamma (\beta) \Psi(\beta)$, and therefore
$$
f'_1(\beta, t)=-\frac{\ln t- \Psi (2-\beta)}{\lambda t^\beta
    \Gamma (2-\beta)}.
$$
In order to dispense with the singularity in the denominators of functions $f_1(\beta, t)$ and $f'_1(\beta, t)$, we use the relations

\[
\frac{1}{\Gamma(2-\beta)}=\frac{2-\beta}{\Gamma(3-\beta)}, \quad
\Psi(2-\beta)=\Psi(3-\beta)-\frac{1}{2-\beta}.
\]
Then
\begin{equation}\label{ash33}
f_1(\beta, t)=\frac{(2-\beta)}{\lambda t^\beta\Gamma (3-\beta)},
\end{equation}
and the function  $f'_1(\beta, t)$ can be rewritten as
\begin{equation}\label{ash34}
f'_1(\beta, t)=-\frac{1}{\lambda t^\beta} \cdot\frac{(2-\beta)[\ln
    t - \Psi (3-\beta)]-1}{ \Gamma (3-\beta)}.
\end{equation}

If $\gamma\approx 0,57722$ - is the Euler-Masheroni constant, then
$-\gamma <\Psi(3-\beta)< 1-\gamma$. Therefore, for sufficiently large $t$ and  $\beta\in [\beta_1, \beta_2]$ the following estimate is true
\begin{equation}\label{ash35}
-f'_1(\beta, t)\geq \frac{1}{\lambda t^\beta}.
\end{equation}

We pass to the estimate of derivative  $f'_2(\beta, t)$. We denote the integrand in  (\ref{ash32})  through $F(\zeta, \rho)$:
\[
F(\zeta, \beta)=\frac{1}{2\pi i \beta\lambda t^\beta}\cdot
\frac{e^{\zeta^{1/\beta}}\zeta^{1-1/\beta}}{\zeta+\lambda
    t^\beta}.
\]
Note that the domain of integration  $\delta(\theta)$ also depends on  $\beta$. To consider this when differentiating function
$f'_2(\beta, t)$, we rewrite the integral (\ref{ash32}) in the following form:
\[
f_2(\beta, t)=f_{2+}(\beta, t)+f_{2-}(\beta, t)+f_{21}(\beta, t),
\]
where
\[
f_{2\pm}(\beta, t)=e^{\pm i \theta}\int\limits_1^\infty
F(s\,e^{\pm i \theta}, \beta)\, ds,
\]

\[
f_{21}(\beta, t) = i \int\limits_{-\theta}^{\theta} F(e^{i y},
\beta)\, e^{iy} dy= i\theta \int\limits_{-1}^{1} F(e^{i \theta s},
\beta)\, e^{i\theta s} ds.
\]

Let us consider function   $f_{2+}(\beta, t)$. Since $\theta =
\theta_\epsilon \beta$ and $\zeta= s\, e^{i\theta}$,  then
\[
e^{\zeta^{1/\beta}}=e^{-s^{\frac{1}{\beta}}(\varepsilon_1-i\varepsilon_2)},
\quad \cos \theta_\epsilon =-\varepsilon_1<0, \quad \sin
\theta_\epsilon= \varepsilon_2>0.
\]
therefore,
\[
f_{2+}(\beta, t)=\frac{1}{I}\int\limits_1^\infty \frac{J(s,
    \beta)}{e^{i\theta_\epsilon \beta} s +\lambda t^\beta} ds, \quad
I=2\pi i \beta\lambda t^\beta, \quad J(s,
\beta)=e^{-s^{\frac{1}{\beta}}(\varepsilon_1-i\varepsilon_2)}
s^{1-\frac{1}{\beta}} e^{i\theta_\epsilon (2\beta-1)}.
\]
Simple calculations show the following
$$
f_{2+}'(\beta, t)=\frac{1}{I}\int\limits_1^\infty \frac{J(s,
    \beta)\big[-\frac{\varepsilon_1-i\varepsilon_2}{\beta^2}s^{1/\beta}\ln
    s+\frac{\ln s}{\beta^2}+2i\theta_\epsilon -\frac{1}{\beta}-\ln
    t-\frac{i\theta_\epsilon e^{i\theta_\epsilon \beta}  s\,+\lambda
        t^\beta \ln t}{e^{i\theta_\epsilon \beta} s +\lambda
        t^\beta}\big]}{e^{i\theta_\epsilon \beta} s +\lambda t^\beta} ds.
$$
We have  $|e^{i\theta_\epsilon \beta} s +\lambda t^\beta|\geq
\varepsilon_3\lambda t^\beta$, $\sin (\pi-\theta_\epsilon)=\sin
\epsilon =\varepsilon_3>0$. Therefore
$$
|f_{2+}'(\beta, t)|\leq \frac{1}{\varepsilon_3|I|\lambda
    t^\beta}\int\limits_1^\infty
e^{-\varepsilon_1\,s^{\frac{1}{\beta}}}s^{1-\frac{1}{\beta}}\,\bigg[\frac{1}{\beta^2}(s^{1/\beta}+1)\ln
s+2\theta_\epsilon +\frac{1}{\beta}+(1+\frac{1}{\varepsilon_3})\ln
t+\frac{\theta_\epsilon s}{\varepsilon_3\lambda t^\beta}\bigg]
ds
$$
$$
\leq\frac{C}{(\varepsilon_3 \lambda t^\beta)^2}\int\limits_1^\infty
e^{-\varepsilon_1\,s^{\frac{1}{\beta}}}s^{1-\frac{1}{\beta}}\,(s+\ln
t)\, ds\leq \frac{C}{(\varepsilon_3 \lambda
    t^\beta)^2}\int\limits_1^\infty
e^{-\varepsilon_1\,s^{\frac{1}{2}}}s^{\frac{1}{2}}\,(s+\ln t)\,
ds\leq \frac{C(\varepsilon_1^{-5} + \varepsilon_1^{-3}\ln
    t)}{(\varepsilon_3 \lambda t^\beta)^2}.
$$
Derivative  $f'_{2-}(\beta, t)$ has a similar estimate.

Let us estimate the derivative of function  $f_{21}(\beta, t)$. We have
\[
f'_{21}(\beta, t)=\frac{\theta_\epsilon}{2\pi \lambda
    t^\beta}\int\limits_{-1}^1\frac{e^{e^{i\theta_\epsilon
            s}}e^{i\theta_\epsilon s(\beta-1)}\big[i\theta_\epsilon s-\ln
    t-\frac{i\theta_\epsilon e^{i\theta_\epsilon \beta} s+\lambda
        t^\beta \ln t}{e^{i\theta_\epsilon \beta} s +\lambda
        t^\beta}\big]}{e^{i\theta_\epsilon \beta} s +\lambda t^\beta} ds.
\]
Therefore,
\[
|f'_{21}(\beta, t)|\leq\frac{C(1+ \ln t)}{(\lambda t^\beta)^2}.
\]

So, from the estimates of functions  $f'_{2\pm}(\beta, t)$ and
$f'_{21}(\beta, t)$, we obtain
\begin{equation}\label{ash36}
|f'_{2}(\beta, t)|\leq C(\beta_2)\cdot\frac{\ln t}{(\lambda
    t^\beta)^2}.
\end{equation}

Comparing estimates  (\ref{ash35}) and (\ref{ash36}) we conclude that there exists a constant  $T_2=T_2(\lambda_1, \beta_2)$ such that for $t\geq T_2$ derivative $\frac{d}{d\beta} e_{\lambda,2}(\beta)$
is negative, that is, function  $e_{\lambda,2}(\beta)$ is decreasing in  $\beta\in [\beta_1, \beta_2]$ for all  $\lambda\geq
\lambda_1$.

To complete the proof of the lemma, note that, using reasoning similar to the above, we prove the estimates
\[
|f_{2\pm}(\beta, t)|\leq\frac{C}{\varepsilon_1^3 \varepsilon_3
    (\lambda t^\beta)^2}, \quad |f_{21}(\beta,
t)|\leq\frac{C}{(\lambda t^\beta)^2},
\]
and therefore
\begin{equation}\label{ash37}
|f_{2}(\beta, t)|\leq C(\beta_2)\cdot\frac{1}{(\lambda
    t^\beta)^2}.
\end{equation}


In the next lemma, in contrast to Lemma  \ref{e2}, we assume that $\beta_1> 1$ and denote $\varepsilon_4=\beta-1$
.

\textbf{Lemma.}\label{e1}
    Let $\beta_1> 1$ and $\beta_2<2$. There exists a positive number
    $T_1=T_1(\lambda_1, \beta_1, \beta_2)$ such that for all  $t\geq
    T_1$ and $\lambda \geq \lambda_1$ function $e_{\lambda,1}(\beta) =
    E_{\beta, 1}\big(-\lambda t^\beta\big)$ is negative and monotonically increasing in  $\beta\in [\beta_1, \beta_2]$  and the following estimates are true
    \[
    e_{\lambda,1}(\beta) =-\frac{1}{\lambda
        t^\beta}\cdot\frac{(\beta-1)(2-\beta)}{\Gamma
        (3-\beta)}+\frac{O(1)}{(\lambda t^\beta)^2},
    \]
    \[
    \frac{d}{d\beta} e_{\lambda,1}(\beta) =\frac{1}{\lambda
        t^\beta}\cdot \frac{(\beta-1)(2-\beta)[\ln t - \Psi
        (3-\beta)]+2\beta-3}{ \Gamma (3-\beta)}+\frac{O(1)}{(\lambda
        t^\beta)^2},
    \]
    where $O(1)$ depends only on $\beta_2$.


\textbf{Proof.} Let the numbers   $\delta(\theta)$, $\epsilon$,
    $\theta_\epsilon$ and $\varepsilon_j$ be defined in the same way as above. The following formula is true (see \cite{as6}, p. 135)
    \begin{equation}\label{ash38}
    E_{\beta, 1}\big(-\lambda t^\beta\big)= \frac{1}{\lambda t^\beta \Gamma
        (1-\beta)}-\frac{1}{2\pi i \beta\lambda
        t^\beta}\int\limits_{\delta(\theta)}\frac{e^{\zeta^{1/\beta}}\zeta}{\zeta+\lambda
        t^\beta} d\zeta = g_1(\beta, t)+g_2(\beta, t).
    \end{equation}

    Using the above properties of functions  $\Gamma$ and $ \Psi$,
    we rewrite function  $g_1(\beta, t)$ in the following form
    \[
    g_1(\beta, t)=-\frac{(\beta-1)(2-\beta)}{\lambda t^\beta \Gamma
        (3-\beta)},
    \]
    and its derivative
    $$
    g'_1(\beta, t)=-\frac{\ln t- \Psi (1-\beta)}{\lambda t^\beta
        \Gamma (1-\beta)},
    $$
    in the form
    $$
    g'_1(\beta, t)=\frac{(\beta-1)(2-\beta)[\ln t- \Psi (3-\beta)]
        +2\beta-3}{\lambda t^\beta \Gamma (3-\beta)}.
    $$
    It is clear that number  $T_1$ can be chosen so large, depending on $\varepsilon_4$, that for $t\geq T_1$ the following estimate is true
    \begin{equation}\label{ash39}
    g'_1(\beta, t)\geq \frac{1}{\lambda t^\beta }.
    \end{equation}

Defining function $G(\zeta, \beta)$   by the following equation
\begin{equation}\label{ash40}
G(\zeta, \beta)=\frac{1}{2\pi i \beta\lambda t^\beta}\cdot
\frac{e^{\zeta^{1/\beta}}\zeta}{\zeta+\lambda t^\beta},
\end{equation}
introduce functions  $g_{2\pm}(\beta, t)$ and $g_{21}(\beta, t)$ similar to the above. Then we have
$$
g_{2}(\beta, t)=g_{2+}(\beta, t)+g_{2-}(\beta, t)+g_{21}(\beta,
t).
$$
Repeating the same reasoning as above, we get
$$
|g_{2\pm}'(\beta, t)|\leq\frac{C}{(\varepsilon_3 \lambda
    t^\beta)^2}\int\limits_1^\infty
e^{-\varepsilon_1\,s^{\frac{1}{\beta}}}s\,(s+\ln t)\, ds\leq
\frac{C(\varepsilon_1^{-6} + \varepsilon_1^{-4}\ln
    t)}{(\varepsilon_3 \lambda t^\beta)^2}.
$$
Likewise,
\[
|g'_{21}(\beta, t)|\leq\frac{C(1+ \ln t)}{(\lambda t^\beta)^2}.
\]
Therefore, from the estimates of functions  $g'_{2\pm}(\beta, t)$ and
$g'_{21}(\beta, t)$, we obtain
\begin{equation}\label{ash41}
|g'_{2}(\beta, t)|\leq C(\beta_2)\cdot\frac{\ln t}{(\lambda
    t^\beta)^2}.
\end{equation}

Comparing estimates  (\ref{ash39}) and (\ref{ash41}) we conclude that there exists a constant  $T_1=T_1(\lambda_1, \beta_1, \beta_2)$ such that for  $t\geq T_1$ the derivative of  $\frac{d}{d\beta}
e_{\lambda,1}(\beta)$ is positive, that is, function
$e_{\lambda,1}(\beta)$ is increasing in   $\beta\in
[\beta_1, \beta_2]$  for all $\lambda\geq \lambda_1$.

It is easy to check the validity of the estimates
\[
|g_{2\pm}(\beta, t)|\leq\frac{C}{\varepsilon_1^4 \varepsilon_3
    (\lambda t^\beta)^2}, \quad |g_{21}(\beta,
t)|\leq\frac{C}{(\lambda t^\beta)^2},
\]
and therefore,
\begin{equation}\label{ash42}
|g_{2}(\beta, t)|\leq C(\beta_2)\cdot\frac{1}{(\lambda
    t^\beta)^2}.
\end{equation}


The following lemma, in a slightly different formulation, was proven in \cite{as27}. For the convenience of the readers, here are the main points of the proof.

\textbf{Lemma.}\label{e3}
    Let $\alpha_1> 0$. There exists a positive number
    $T_3=T_3(\lambda_1, \alpha_1)$ such that for all $t\geq T_3$ and
    $\lambda \geq \lambda_1$ function  $e_{\lambda,3}(\alpha) =
    E_{\alpha, 1}(-\lambda t^\alpha)$ is positive and monotonically decreasing in  $\alpha\in [\alpha_1, 1)$ and the following estimates are true
    \[
    e_{\lambda,3}(\alpha) =\frac{1}{\lambda
        t^\alpha}\cdot\frac{1-\alpha}{\Gamma
        (2-\alpha)}+\frac{O(1)}{(\lambda t^\alpha)^2},
    \]
    \[
    \frac{d}{d\alpha} e_{\lambda,3}(\alpha) = - \frac{1}{\lambda
        t^\alpha}\cdot \frac{(1-\alpha)[\ln t - \Psi (2-\alpha)]+1}{
        \Gamma (2-\alpha)}+\frac{O(1)}{(\lambda t^\alpha)^2},
    \]
    where $O(1)$ depends only on  $\alpha_1$.

\textbf{Proof.} The parameter of contour
    $\delta(\theta)$ is chosen as $\theta=\frac{3\pi}{4}\alpha$,
    $\alpha\in [\alpha_1, 1)$. Then, by the definition of a contour, we have (see \cite{as6}, p.135)
    \begin{equation}\label{ash43}
    E_{\alpha, 1}(-\lambda t^\alpha)= \frac{1}{\lambda t^\alpha \Gamma
        (1-\alpha)}-\frac{1}{2\pi i \alpha\lambda
        t^\alpha}\int\limits_{\delta(\theta)}\frac{e^{\zeta^{1/\alpha}}\zeta}{\zeta+\lambda
        t^\alpha} d\zeta = q_1(\alpha, t)+q_2(\alpha, t).
    \end{equation}

Since function  $\Gamma(1-\alpha)$ tends to infinity for  $\alpha\rightarrow 1$
it is convenient to represent function  $q_1(\alpha, t)$ in the following form

$$
q_1(\alpha, t)=\frac{1-\alpha}{\lambda t^\alpha \Gamma
    (2-\alpha)}.
$$

The derivative of function  $q_1(\alpha, t)$ has the form
$$
q_1'(\alpha, t)=-\frac{\ln t - \Psi (1-\alpha)}{\lambda t^\alpha
    \Gamma (1-\alpha)}.
$$
Therefore ,
\begin{equation}\label{ash44}
-q_1'(\alpha, t)=\frac{1}{\lambda t^\alpha}\cdot
\frac{(1-\alpha)[\ln t - \Psi (2-\alpha)]+1}{ \Gamma
    (2-\alpha)}\geq \frac{1}{\lambda t_0^{\alpha}}.
\end{equation}

Note that the integrand in equation (\ref{ash43}) is $G(\zeta, \alpha)$ (see (\ref{ash40})). By virtue of the choice of $\theta$, along the contour $\delta(\theta)$, the following relations hold
\[
\zeta= s\, e^{i\theta}, \quad e^{\zeta^{1/\alpha}}=e^{\frac{1}{2}
    (i-1)\,s^{\frac{1}{\alpha}}},\quad |\zeta+\lambda t^\alpha|\geq
\lambda t^\alpha.
\]
Therefore, repeating similar reasoning as in the proof of Lemma  \ref{e1} (see. \cite{as27})), we obtain the following estimates
\[
|q_2(\alpha, t)|\leq \frac{C}{(\lambda t^{\alpha})^2},
\]
\[
|q_2'(\alpha, t)|\leq C \,\frac{\ln t+\alpha_1^{-1}}{(\lambda
    t^{\alpha})^2}\leq C(\alpha_1)\, \frac{\ln t}{(\lambda
    t^{\alpha})^2}.
\]

Comparing the last estimate with the estimate in  (\ref{ash44}), we conclude that there exists a constant  $T_3=T_3(\lambda_1, \alpha_1)$ al) such that for
$t\geq T_3$ derivative $\frac{d}{d\alpha} e_{\lambda,3}(\alpha)$
is negative, that is, function  $e_{\lambda,3}(\alpha)$ is decreasing in  $\alpha\in [\alpha_1, 1)$ for all  $\lambda\geq\lambda_1$.

In conclusion, note that the inequality  $e_\lambda(1)=
e^{-\lambda t} >0$ implies that  $e_\lambda (\alpha)$
is positive for all   $\alpha \in [\alpha_1, 1)$.


\section{INVERSE PROBLEM}

The proof of Theorem \ref{it} is based on the monotonicity properties of the Mittag - Leffler functions $E_{\rho, \mu} (t)$ in parameter $\rho$, proven in the previous section.

Note that conditions (\ref{ash12}) and (\ref{ash13}) can be considered as equations to determine parameters $\alpha$ and $\beta$. Moreover, as noted above, equation (\ref{ash13}) includes only unknown  one $\beta$ , and equation (\ref{ash12}) includes two unknowns $\alpha$ and $\beta$ (see the solution forms in (\ref{ash10}) and (\ref{ash11})). Therefore, we begin by considering condition (\ref{ash13}), which is convenient to write in the form of (\ref{ash14}) (note that $\varphi_{k_0}\neq 0$))

\[ P(\beta)
\varphi_{k_0}=d_2, \quad P(\beta)\equiv\frac{\Delta_{k_0}(t_2,
    \beta)}{\Delta_{k_0}(T, \beta)}.
\]
In order to prove the existence of a unique parameter $\beta$, satisfying this equation, it is sufficient to show that the derivative of function $P(\beta)$ maintains its sign.

Let us introduce the notation for the pair of functions $p(\beta)$ and $q(\beta)$:
$V(p,q)\equiv p'q-pq'$. Then
\[
P'(\beta) =\frac{V (\Delta_{k_0}(t_2, \beta),\,\, \Delta_{k_0}(T,
    \beta))}{\Delta^2_{k_0}(T, \beta)}.
\]
Now it suffices to show that  $V (\Delta_{k_0}(t_2,\,\,
\beta), \Delta_{k_0}(T, \beta))$ maintain the sign. As required in Theorem  \ref{it}, we consider only sufficiently large values of $t_2$. By virtue of the asymptotic estimates established in Lemmas \ref{e2} and \ref{e1}, the principal part of the function $\Delta_{k_0}(t,
\beta)$ for
$$
t\geq T_{1,2}=T_{1,2} (\lambda_1, \beta_1, \beta_2)\equiv\max \{\,
T_{1} (\lambda_1, \beta_1, \beta_2), \quad T_{2} (\lambda_1,
\beta_2)\,\}
$$
has the following form
\[
\lambda_{k_0} t f_1(\beta, t) - g_1(\beta, t).
\]
 As follows from the definitions of these functions, due to the presence of factor $\lambda_{k_0} t$, to prove the unique solvability of equation (\ref{ash14}), it is sufficient to check the sign of function $V \big(\lambda_{k_0} t_2 f_1(t_2, \beta),\,\, \lambda_{k_0} T
 f_1(T, \beta)\big)$. It is easy to check that
\[
V \big(\lambda_{k_0} t_2 f_1(t_2, \beta),\,\, \lambda_{k_0} T f_1(T,
\beta)\big)=\frac{1}{(t_2 T)^{\beta-1}}\cdot
\frac{2-\beta}{\Gamma(3-\beta)}\cdot \Big[\ln T-\ln t_2
+\frac{O(1)}{\ln t_2}\Big],
\]
where $O(1)$ depends only on $\beta_2$.  Therefore, for $T_{1,2} \leq t_2< T $ function $P(\beta)$ has a strictly positive derivative, which, in turn, means that this function is strictly monotonic, in particular, $P(\beta_1) \leq P(\beta)
\leq P(\beta_2)$, $\beta\in [\beta_1, \beta_2]$. This implies the unique solvability of equation (\ref{ash14}) with respect to $\beta$.

Note that the above proof implies that condition (\ref{ash15}) for a given number $d_2$ is sufficient for the unique solvability of equation (\ref{ash14}). Thus, it is proven that for sufficiently large $t_2$ this condition is not only necessary but also sufficient.

Now let $t_1 \geq T_3$ (see Lemma \ref{e3}). Then the unique solvability of equation (\ref{ash14})  with respect to $\alpha$, under condition (\ref{ash16}), follows directly from Lemma \ref{e3} and the following equation

\[
W(\alpha, \beta)=\int\limits_\Omega |u(x,t_1)|^2
dx=\sum\limits_{k=1}^\infty E_{\alpha, 1}(-\lambda_k
t^\alpha)\bigg|\frac{ \varphi_k}{\Delta_k}\bigg|^2.
\]
Thus, the proof of Theorem \ref{it} is complete.

The authors are deeply grateful to Sh.A. Alimov for discussing the results of the research.


%
%

\end{document}